\documentclass[11pt]{amsart}
\usepackage{}
\usepackage{amssymb}
\theoremstyle{plain}
\newtheorem{Thm}{Theorem}

\newtheorem{Pro}[Thm]{Proposition}

\begin{document}

\title[Expanding Ricci solitons with pinched curvature]
{Expanding Ricci solitons with pinched Ricci curvature}

\author{Li Ma }

\address{Li Ma, Department of Mathematical Sciences, Tsinghua University,
 Peking 100084, P. R. China}

\email{lma@math.tsinghua.edu.cn}

\thanks{The research is partially supported by the National Natural Science
Foundation of China 10631020 and SRFDP 20090002110019}

\begin{abstract}
In this paper, we prove that expanding gradient Ricci solitons
with (positively) pinched Ricci curvature are trivial ones.
Namely, they are either compact or flat.


\end{abstract}

 \maketitle

\section{Introduction}
In this paper we consider Problem 9.62 in the famous book
\cite{Cho06}, which is also the unanswered question in \cite{MD}.
Namely, when $n\geq 3$, do there exist expanding gradient Ricci
solitons with (positively) pinched Ricci curvature? In dimension
three, we settle the problem completely. Here we recall that the
(positively) pinched Ricci curvature for the Riemannian manifold
$(M^n,g)$ is in the sense that
\begin{equation}\label{pinch}
Rc\geq \epsilon Rg\geq 0, \end{equation} where $R$ and $Rc$ are
the scalar and Ricci curvatures of the metric $g$ respectively,
$\epsilon>0$ is a small constant. This concept plays an important
role in the seminal work of R.Hamilton \cite{Ham82}. We remark
that the compact expanding gradient Ricci solitons are Einstein.
This result is known in G.Perelman \cite{P02}. Then we may assume
that the expanding gradient Ricci soliton $(M^n,g(t),\phi)$ under
consideration is complete, non-compact, and in canonical form that
$$
D^2\phi=Rc+\frac{1}{2t}g, \ \ t>0.
$$
We denote by $d_g(x,o)$  the distance between the points $x$ and
$o$ in $(M,g)$.

We show that there are only trivial ones.
\begin{Thm}\label{thm:1}
Expanding gradient Ricci solitons (M,g) with (positively) pinched
Ricci curvature and curvature decay at the order
$d_g(x,o)^{-2-\delta}$ for some $\delta>0$, are trivial ones.
Namely, they are either compact or $R^n$ with flat metric.
\end{Thm}

We remark that in dimension three, the curvature decay condition
is automatically true \cite{MD}. This result is used in \cite{M}.
In the dimension bigger than three, the same is true for locally
conformally flat expanding gradient Ricci solitons with pinched
Ricci curvature. In general, since we are studying the Ricci flow,
we should have the curvature decay order as that of the Ricci
curvature. This will be considered in the future.

This paper is organized as follows. In section \ref{sect2} we
recall some famous results, which will be in use in section
\ref{sect3}. Theorem \ref{thm:1} is proved in section \ref{sect3}.

We shall use $r$ denote various uniform positive constants.

\section{Preliminary}\label{sect2}
Before we prove our main result Theorem \ref{thm:1}, we cite the
following results, which will be in use in next section. The first
is

\begin{Pro}\label{pro:1}(Hamilton, Proposition 9.46 in
\cite{Cho06}).
If $(M,g(t))$, $t>0$, is a complete non-compact expanding gradient
Ricci soliton with $Rc>0$, then
$$
AVR(g(t)):=\lim_{r\to \infty}\frac{B_{g(t)}(o,r)}{r^n}>0,
$$
where the definition of $AVR(g(t))$ is independent of the base
point $o\in M$.
\end{Pro}

The second is

\begin{Pro}\label{pro:2} (\cite{MD}). If $(M,g(t))$, $t>0$, is a complete non-compact expanding gradient
Ricci soliton with (\ref{pinch}), then the scalar curvature is
quadratic exponential decay.
\end{Pro}

The third one is Theorem 1.1 in \cite{B}. Roughly speaking, the
result says that for the complete and non-compact $(M,g)$, if
$AVR(g)>0$ and the curvature decay suitably, then it is an
asymptotic manifold. We invite the readers to papers \cite{B} and
\cite{LP} for the definitions of asymptotic flat manifold
$M_{\tau}$ and coordinates $(z)=(z^i)$ at infinity (also called
the asymptotic coordinates).

The last one is Proposition 10.2 in \cite{LP}. Namely,
\begin{Pro}\label{pro:3} (\cite{LP}).If $(M,g)$ is asymptotic flat
with $g\in \mathbf{M_{\tau}}$, for some $\tau>\frac{n-2}{2}$, and
the Ricci curvature is non-negative. Then the mass
$$
m(g):=\lim_{r\to\infty}\int_{S_r}\mu\lrcorner dz
$$
 is non-negative, with $m(g)=0$
if and only if $(M,g)$ is isometric to $R^n$ with its Euclidean
metric. Here $S_r:\{x\in M; d_g(o,x)=1\}$,  $$
\mu=(\partial_ig_{ij}-\partial_jg_{ii})\partial_j, \ \
\partial_j=\frac{\partial}{\partial z^j},
$$
and $(z^j)$ are the asymptotic coordinates.
\end{Pro}

\section{Proof of Theorem \ref{thm:1}}\label{sect3}
We argue by contradiction. So we assume that $(M,g(t))$ is not
flat.

Using the strong maximum principle \cite{Shi89a} to the Ricci
soliton, we may assume that the scalar curvature is positive,
i.e., $R>0$. Hence we know that $Rc>0$. According to the arguments
in \cite{Cho06} and \cite{MD}, we know that $\phi$ is a proper
strict convex function, which implies by using the Morse theory
that $M^n$ is diffeomorphic to $R$. Using Proposition \ref{pro:1},
Proposition \ref{pro:2}, and Theorem 1.1 in \cite{B} we know that
$(M,g(t))$ is an asymptotic flat manifold. We also know that
$$
\phi(x)\approx d_g(x,o)^2, \ \ |\nabla \phi(x)|\approx d_g(x,o).
$$

Recall that in coordinates $(z^j)$ at infinity, we have the Ricci
soliton equation
$$
\phi_{ij}=R_{ij}+\frac{1}{2t}g_{ij}.
$$
For notation simple, we let $t=1/2$. Then we have
$g_{ij}=\phi_{ij}-R_{ij}$. By Ricci pinching condition we know
that $R_{ij}$ decay exponentially and $$\Delta
\phi=\phi_{ii}=R+n.$$ Recall the Ricci formula
$$
\phi_{iji}=\phi_{iij}+R_{ji}\phi_i.
$$
This also implies that
$$
\nabla R=-2Rc(\nabla \phi).
$$
Hence $|\nabla R|$ decays in the exponent rate.

 We now
compute the mass. Using the Ricci formula, we have
\begin{eqnarray*}
m(g)&=&\lim_{r\to\infty}\int_{S_r}(\partial_ig_{ij}-\partial_jg_{ii})\partial_j\lrcorner
dz
\\
&=&\lim_{r\to\infty}\int_{S_r}(\partial_i\phi_{ij}-\partial_j\phi_{ii})\partial_j\lrcorner
dz\\
&=&\lim_{r\to\infty}\int_{S_r}\phi_{iji}\partial_j\lrcorner dz\\
&=&\lim_{r\to\infty}\int_{S_r}(\phi_{iij}+R_{ij}\phi_i)\partial_j\lrcorner dz\\
&=&\lim_{r\to\infty}\int_{S_r}R_j\partial_j\lrcorner dz\\
&=&0.
\end{eqnarray*}
Using Proposition \ref{pro:3} we know that $(M,g(1/2))$ is $R^n$
with the Euclidean metric. A contradiction.
 This completes the
proof of Theorem \ref{thm:1}. q.e.d.


\begin{thebibliography}{20}

\bibitem{B}
S.Bando, A.Kasue, H.Nakajima,\emph{On a construction of
coordinates at infinity on manifolds with fast curvature decay and
maximal volume growth}, Invent. math. 97, 313-349 (1989).

\bibitem{Cho06}
 B.Chow, P.Lu, L.Ni, \emph{Hamilton's Ricci Flow}. Science
Press. American Mathematical Society, Beijing,Providence, 2006.



\bibitem{Ham82}
R.Hamilton, \emph{Three-manifolds with positive Ricci curvature}.
 J. Differential Geom., 2(1982)255-306.



\bibitem{LP}
J.Lee, T.Parker, \emph{The Yamabe problem}, Bull. A.M.S.,
17(1)(1987)37-91.

\bibitem{M}
Li Ma, \emph{A proof of Hamilton's conjecture}. arXiv:1003.0177



\bibitem{MD}
Li Ma, DeZhong Chen, \emph{Remarks on complete non-compact
gradient Ricci expanding solitons}, Kodai Math.Journal, 2010,
online



\bibitem{P02} Grisha Perelman,
\emph{The entropy formula for the Ricci flow and its geometric
applications}, http://arxiv.org/abs/math/0211159v1


\bibitem{Shi89a}
W. X. Shi, \emph{Deforming the metric on complete Riemannian
manifolds},JDG, 30(1989)223-301.



\end{thebibliography}
\end{document}